\documentclass[12pt]{amsart}
\usepackage[a4paper]{geometry}
 \newtheorem{thm}{Theorem}[section]
 
 \newtheorem{lem}[thm]{Lemma}
 \newtheorem{prop}[thm]{Proposition}
 \theoremstyle{definition}
 \newtheorem{defn}{Definition}[section]
 \theoremstyle{remark}
 
 \numberwithin{equation}{section}

\begin{document}
\title {\textbf{Biharmonic $\delta(\lowercase{r})$-ideal hypersurfaces in Euclidean spaces are minimal}}
\author{Deepika and Andreas Arvanitoyeorgos}

 \maketitle
% -----------------------------------------------------------
\begin{abstract}
A submanifold $M^n$ of a Euclidean space $\mathbb{E}^N$ is called
biharmonic if $\Delta\vec{H}=0$, where $\vec{H}$ is the mean
curvature vector of $M^n$.  A well known conjecture of B.Y. Chen
states that the only biharmonic submanifolds of Euclidean spaces are
the minimal ones. Ideal submanifolds were introduced by Chen as
those which receive the least possible tension at each point. In
this paper we prove that every $\delta(r)$-ideal biharmonic
hypersurfaces in the Euclidean space $\mathbb{E}^{n+1}$ ($n\geq 3$)
is minimal. In this way we generalize a recent result of B. Y. Chen
and M. I. Munteanu. In particular, we show that every
$\delta(r)$-ideal biconservative hypersurface in Euclidean space
$\mathbb{E}^{n+1}$ for $n\geq 3$ must be of constant mean curvature.
\\
\\
\textbf{AMS 2010 Subject Classification:} 53D12, 53C40, 53C42\\
\textbf{Key Words:}  biharmonic submanifolds; biharmonic map;
$\delta$-invariant; $\delta(r)$-ideal submanifolds;
 biconservative hypersurface; mean curvature vector.
\end{abstract}

%%% ----------------------------------------------------------------------
\maketitle
%%% ----------------------------------------------------------------------
\section{\textbf{Introduction}}

In the middle of 1980's B.Y. Chen initiated the investigation of
properties of submanifolds in Euclidean spaces $\mathbb{E}^N$. For
an historical overview we refer to his books \cite{BY1} and
\cite{BY2}. Among several important problems that Chen had raised
that time was the following well-known conjecture \cite{BY7}:

\noindent
 \emph{The only biharmonic submanifolds of Euclidean spaces are the minimal submanifolds.}

 A biharmonic submanifold $M$ is defined by the equation $\Delta\vec{H}=0$, where $\Delta$ and $\vec{H}$ denote respectively the Laplace operator and the mean curvature vector field of $M$.
 It is well known that $M$ is biharmonic if and only if the immersion $x: (M, g)\to\mathbb{E}^N$ is a biharmonic map.

The conjecture was originally proved for surfaces in $\mathbb{E}^3$
by B.Y. Chen in \cite{BY7} and for certain submanifolds in
$\mathbb{E}^n$ (including one dimensional) by I. Dimitri\'c in
\cite{ID}. An alternative approach was proposed by  T. Hasanis and
T. Vlachos in \cite{TT} who proved the conjecture for hypersurfaces
in Euclidean 4-spaces.  Since then several researchers have made
significant contributions towards proving it, such as F. Defever
\cite{DE}, K. Akatagawa and S. Maeta \cite{AM}, Y. Fu \cite{YF2},
\cite{YF3}, R. Shankar and A. Sharfuddin \cite{RA}, and more
recently in B.Y. Chen \cite{BY9}, and N. Koiso, H. Urakawa
\cite{KU}.

In contrast to its simple statement, the conjecture has turned quite endure to several attempts for its proof.  Therefore, it was natural for researchers to impose some natural assumptions.
The most usual one, and in fact the most successful into confirming the conjecture for several cases, was to assume that the submanifold is a hypersurface in the Euclidean space.  In this case one usually makes some extra assumption about the number of distinct eigenvalues of the shape operator, or about  the scalar curvature.

Another natural assumption for the hypersurface $M$ is to be an {\it ideal} (or {\it $\delta(r)$-ideal}) hypersurface.
We give the formal definition in Section 2.  Such hypersurfaces were introduced by B.Y. Chen via the concept of $\delta$-invariants, in his investigation to define ``nice immersed submanifolds" as those which receive the least possible tension at each point.
We refer to \cite{BY2} and \cite{BY8} for a deeper motivation.

In the work \cite{BY3} B.Y. Chen and M. I. Munteanu proved that $\delta(2)$-ideal and $\delta(3)$-ideal biharmonic hypersurfaces of a Euclidean space is minimal.
In the present work we extend this result by proving the following:

\begin{thm}\label{main1}
Every $\delta(r)$-ideal oriented biharmonic hypersurface with at most $r+1$ distinct principal curvatures in the Euclidean
spaces $\mathbb{E}^{n+1}$ ($n\geq 3$), is minimal.
\end{thm}

Closely related to the concept of biharmonic submanifolds is the concept of {\it biconservative} submanifolds.
These were introduced by R. Caddeo et al. in \cite{RSC} and are submanifolds  with conservative stress-energy tensor.

For the case of hypersurfaces $M^n$ in $\mathbb{E}^{n+1}$ it can be shown that the biconservativity condition
 is equivalent to the equation
$2\mathcal{A}({\rm grad}H)+nH{\rm grad}H=0$, where $\mathcal{A}$ is the shape operator of $M^n$ and $H$ the mean curvature.
This equation is one of the two equations which are equivalent to the condition of biharmonicity,
$\Delta\vec{H}=0$ (cf. (\ref{1.2}), (\ref{1.1})).
Therefore, a biharmonic  hypersurface is biconservative.

Biconservative hypersurfaces had appeared in the literature under the name {\it H-hypersurfaces} (\cite{TT}).
They have attracted recently the interest of several researchers (e.g.
\cite{DEEP}, \cite{DEEP1},  \cite{DCA},  \cite{YF}, \cite{RSG},
\cite{FNA},  \cite{SCA}, \cite{SCA2},  \cite{NC}).
From the proof of Theorem \ref{main1} we also obtain  the following result:

  \begin{prop}\label{main}
Every $\delta(r)$-ideal oriented biconservative hypersurface with at most $r+1$ distinct principal curvatures in Euclidean
spaces $\mathbb{E}^{n+1}$ ($n\geq 3$), has constant mean curvature.
\end{prop}

 We briefly present the central idea of the proof of Theorem \ref{main1} which is
simpler than the method used in \cite{BY3} and \cite{DEEP1}.
Using that $M^n$ is a $\delta(r)$ ideal hypersurface its shape operator has a simpler form.
Since $M^n$ is biharmonic in particular it is biconservative, hence we use the corresponding equation
to simplify the  connection forms by using Codazzi
equation and Gauss equation. Then we see that the definition of mean
curvature provides us an equation showing the relation between
eigenvalues of the shape operator and mean curvature $H$. This equation
plays a very important role in the proof. By
differentiating this equation two or more times we obtain
polynomial equations showing relations among the eigenvalues,
connection forms and mean curvature $H$. Then using a standard
argument involving the resultant of two polynomials as defined in
Lemma \ref{lemma2.2}, we are able to eliminate all the eigenvalues as well as the
connection forms one by one, to obtain an algebraic polynomial
equation in $H$ with constant coefficients which implies that $H$
must be constant.  By taking into account  the second equation that comes from the biharmonicity assumption, we prove that $H$ is zero.

\section{\textbf{Preliminaries}}

   Let ($M^{n}, g$) be an oriented  hypersurface isometrically immersed in Euclidean space
$(\mathbb{E}^{n+1}, \overline g)$, that is $g$ is the induced metric
by the immersion that defines the hypersurface.
  Let $\overline\nabla $ and $\nabla$ denote the linear connections on $\mathbb{E}^{n+1}$ and $M$ respectively. Then the Gauss and Weingarten formulae are given by
\begin{equation}\label{2.1}
\overline\nabla_{X}Y = \nabla_{X}Y + h(X, Y), \hspace{.2 cm} \forall
\hspace{.2 cm}X, Y \in\Gamma(TM),
\end{equation}
\begin{equation}\label{2.2}
\overline\nabla_{X}\xi = -\mathcal{A}_{\xi}X,
\end{equation}
where $\xi$ be the unit normal vector to $M$, $h$ is the second
fundamental form and $\mathcal{A}$ is the shape operator. It is well
known that the second fundamental form $h$ and shape operator
$\mathcal{A}$ are related by
\begin{equation}\label{2.3}
\overline{g}(h(X,Y), \xi) = g(\mathcal{A}_{\xi}X,Y).
\end{equation}

 The mean curvature is given by
\begin{equation}\label{2.4}
H = \frac{1}{n} \mbox{trace} \mathcal{\mathcal{A}},
\end{equation}
and the mean curvature vector $\vec{H}=H\xi$ is a well defined normal vector field to $M^n$ in $\mathbb{E}^{n+1}$.
 The Gauss and Codazzi equations are given by
\begin{equation}\label{2.5}
R(X, Y)Z = g(\mathcal{A}Y, Z) \mathcal{A}X - g(\mathcal{A}X, Z)
\mathcal{A}Y,
\end{equation}
\begin{equation}\label{2.6}
(\nabla_{X}\mathcal{A})Y = (\nabla_{Y}\mathcal{A})X
\end{equation}
respectively, where $R$ is the curvature tensor and
\begin{equation}\label{2.7}
(\nabla_{X}\mathcal{A})Y = \nabla_{X}\mathcal{A}Y-
\mathcal{A}(\nabla_{X}Y)
\end{equation}
for all $ X, Y, Z \in \Gamma(TM)$.\\
The hypersurface $M^n$ is called {\it biharmonic} if
$$
\Delta\vec{H}=0.
$$
By identifying the tangential and normal parts in the above equation, it is known (\cite{BY1}) that it is equivalent to the system

\begin{equation}\label{1.2}
 2\mathcal{A} (\rm grad \emph{H})+ \emph{nH} \hspace{.1 cm} grad \emph{H} =
 0,\end{equation}
 \begin{equation}\label{1.1}
\triangle H +  \emph{H}\hspace{.1 cm}\rm trace (\mathcal{A}^{2}) =
 0,
\end{equation}
where $\triangle$ is the Laplace operator (our sign convention is
such that $\triangle f= - f''$ when $f$ is a function of one real
variable).

\medskip
Next, we recall the concept of $\delta$-invariants and $\delta$-ideal hypersurfaces.  We refer to \cite{BY2} and \cite{BY6}  for more details.

For a Riemannian manifold $M^{n}$ with $n\geq 3$ and an integer $r
\in [2, n-1]$, let $\tau(p)$ be the scalar curvature at $p\in M^n$
and let $\tau(L^{r})$  be the scalar curvature of a linear subspace
$L^{r}$ of dimension $r\geq 3$ of the tangent space $T_{p}(M)$. The
{\it $\delta$-invariant $\delta(r)$ of $M^n$} at $p$ is defined as

\begin{equation}\label{1.3}
\delta(r)(p)=\tau(p)-\inf_{r} \tau(L^{r}).
\end{equation}

For any $n$-dimensional submanifold $M^{n}$ in a Euclidean space
$\mathbb{E}^{m}$ and for an integer $r=2, \dots , n-1$, Chen proved
the following universal sharp inequality

\begin{equation}\label{1.4}
\delta(r)(p)\leq \frac{n^{2}(n-r)}{2(n-r+1)}H^{2}(p),
\end{equation}
where $H^{2}=\langle \vec{H},\vec{H}\rangle$ is the squared mean
curvature.

\begin{defn} A submanifold $M^{n}$ in $\mathbb{E}^{m}$ is called $\delta(r)$-ideal if  equality  in (\ref{1.4}) is satisfied  identically.
\end{defn}

We will need the following result.

\begin{thm}\label{Theorem2.1} {\rm (\cite[Theorem 13.7]{BY2})}  Let $M^{n}$ be a hypersurface in the Euclidean
spaces $\mathbb{E}^{n+1}$. Then for any integer $r=2, \dots , n-1$ it is
\begin{equation}\label{2.9}
\delta(r)\leq \frac{n^{2}(n-r)}{2(n-r+1)}H^{2}.
\end{equation}
 Equality holds at a point $p$ if and only
if there is an orthonormal basis
$\{e_{1},e_{2},\ldots,e_{n}\}$ at $p$ such that the shape
operator is given by
\begin{equation}\label{4.1}
 \mathcal{A}= \left(
                            \begin{array}{ccc}
                               D_r & 0   \\
                              0 & u_{r}I_{n-r} \\
                             \end{array}
                             \right),
\end{equation}
where $D_{r}=$ ${\rm
diag}(\lambda_{1},\lambda_{2},\dots,\lambda_{r})$ and
$u_{r}=\lambda_{1}+\lambda_{2}+\dots+\lambda_{r}$, where $\lambda_{1}, \lambda_{2},\ldots \lambda_{r}$ are the principal curvature functions of $M^n$ at $p$. If this happens at
every point, we call $M^n$ a $\delta(r)$-ideal hypersurface in
$\mathbb{E}^{n+1}$.
 \end{thm}

  Finally, the following algebraic lemma will be useful to our study.

\begin{lem}\label{lemma2.2}{\rm (\cite[Theorem 4.4, pp. 58--59]{KK})} Let D be
a unique factorization domain, and let $f(X) = a_{0}X^{m}
+a_{1}X^{m-1} + \cdots + a_{m}, g(X) = b_{0}X^{n} + b_{1}X^{n-1} +
\cdots + b_{n}$ be two polynomials in $D[X]$. Assume that the
leading coefficients $a_{0}$ and $b_{0}$ of $f(X)$ and $g(X)$ are
not both zero. Then $f(X)$ and $g(X)$ have a non constant common
factor if and only if the resultant $\Re(f, g)$ of $f$ and $g$ is
zero, where
\begin{center}
$\Re(f,g)=
\begin{vmatrix}
  a_{0} & a_{1} & a_{2} & \cdots & a_{m} &   &   &   \\
    & a_{0} & a_{1} & \cdots & \cdots & a_{m} &   &   \\
    &   & \ddots & \ddots & \ddots & \ddots & \ddots &   \\
   &   &   & a_{0} & a_{1} & a_{2} & \cdots & a_{m} \\
  b_{0} & b_{1} & b_{2} & \cdots & b_{n} &   &   &   \\
    & b_{0} & b_{1} & \cdots & \cdots & b_{n} &  &   \\
    &   & \ddots & \ddots & \ddots & \ddots & \ddots &   \\
    &   &   & b_{0} & b_{1} & b_{2} & \cdots & b_{n} \\
\end{vmatrix}.
$
\end{center}
Here there are $n$ rows of $``a"$ entries and $m$ rows of $``b"$ entries.
\end{lem}

\section{\textbf{$\delta(r)$-ideal biharmonic hypersurfaces in $\mathbb{E}^{n+1}$}}

In the present section we will prove Theorem \ref{main1}.

\medskip
\noindent
{\it Proof of Theorem \ref{main1}}.
Let $M^n$ be an oriented  $\delta(r)$-ideal biharmonic
hypersurface in $\mathbb{E}^{n+1} (n > 2)$.
From Theorem
\ref{Theorem2.1} the shape operator (\ref{4.1}) of $M^n$ with respect to
some orthonormal basis $\{e_{1},e_{2},\dots,e_{n}\}$ can be expressed as
\begin{equation}\label{4.2}
\mathcal{A}(e_{i})=\lambda_{i}e_{i}, \quad i=1,2,\dots,n,
\end{equation}
where $\lambda_{i}=\lambda_{1}+\lambda_{2}+\dots+\lambda_{r}$, for
$i=r+1,\dots,n.$ Since we will need to differentiate the principal
curvature functions of $\mathcal{A}$ we need to know that these are
smooth (at least at some connected component).  To this end, we use
an argument given in \cite[Section 3, lines 3-10]{FH}. The set
$M_{\mathcal{A}}$ of all points of $M^n$, at which the number of
distinct eigenvalues of the shape operator $\mathcal{A}$ (i.e. the
principal curvatures) is locally constant, is open and dense in
$M^n$. Therefore, we  can work only on the connected component of
$M_{\mathcal{A}}$ consisting of points where the number of principal
curvatures is at most $r+1$.
 On that connected
component, the principal curvature functions of $\mathcal{A}$ are always smooth.

\smallskip
\noindent
\underline{Claim:}  The mean curvature $H$ of $M^n$ is constant.

 Assume  the contrary  and we will end up into contradiction.
Then there exists  an open connected subset $U$ of $M$
with grad$_{p}H\neq 0$, for all $p\in U$. From (\ref{1.2}) it is easy
to see that grad$H$ is an eigenvector of the shape operator
$\mathcal{A}$ with corresponding principal curvature
$-\frac{nH}{2}$.

 Without lose of generality we choose $e_{1}$ in the direction of
grad$H$, which gives $\lambda_{1}=-\frac{nH}{2}$. We express grad$H$
as
\begin{equation}\label{3.2}
\mbox{grad} H =\sum_{i=1}^{n} e_{i}(H)e_{i}.
\end{equation}
 As we have taken $e_{1}$ parallel to grad$H$, it is
\begin{equation}\label{3}
e_{1}(H)\neq 0, \ \ e_{i}(H)= 0, \hspace{1 cm} i= 2, \dots, n.
\end{equation}
We express
\begin{equation}\label{3.4}
\nabla_{e_{i}}e_{j}=\sum_{k=1}^{n}\omega_{ij}^{k}e_{k}, \hspace{2
cm} i, j = 1, 2, \dots , n.
\end{equation}
Using (\ref{3.4}) and the compatibility conditions
$(\nabla_{e_{k}}g)(e_{i}, e_{i})= 0$, $(\nabla_{e_{k}}g)(e_{i},
e_{j})= 0$, we obtain
\begin{equation}\label{3.5}
\omega_{ki}^{i}=0, \hspace{1 cm} \omega_{ki}^{j}+ \omega_{kj}^{i}
=0,
\end{equation}
for $i \neq j, $ and $i, j, k = 1, 2,\dots, n$.

We set $\lambda_{r+1}=\lambda_{r+2}=\cdots=\lambda_n=\lambda$ and we  consider the following cases:

\smallskip
\noindent \textbf{Case A.} $\lambda_{i}\neq \lambda,\quad i=2,3,
\dots,r.$\\

 Taking $X=e_{i}, Y=e_{j}$, $(i\neq j)$ in (\ref{2.7}) and using (\ref{4.2}),
(\ref{3.4}), we get
$$(\nabla_{e_{i}}\mathcal{A})e_{j}=e_{i}(\lambda_{j})e_{j}+\sum_{k=1}^{n}\omega_{ij}^{k}e_{k}
(\lambda_{j}-\lambda_{k}).$$

Putting the value of $(\nabla_{e_{i}}\mathcal{A})e_{j}$ in
(\ref{2.6}), we find
$$e_{i}(\lambda_{j})e_{j}+\sum_{k=1}^{n}\omega_{ij}^{k}e_{k}
(\lambda_{j}-\lambda_{k})=e_{j}(\lambda_{i})e_{i}+\sum_{k=1}^{n}\omega_{ji}^{k}e_{k}
(\lambda_{i}-\lambda_{k}),$$ whereby taking inner product with
$e_{j}$ and $e_{k}$, we obtain
\begin{equation}\label{3.6}
e_{i}(\lambda_{j})=
(\lambda_{i}-\lambda_{j})\omega_{ji}^{j}=(\lambda_{j}-\lambda_{i})\omega_{jj}^{i},
\end{equation}
\begin{equation}\label{3.7}
(\lambda_{j}-\lambda_{k})\omega_{ij}^{k}=
(\lambda_{i}-\lambda_{k})\omega_{ji}^{k},
\end{equation}
respectively, for distinct $i, j, k = 1, 2, \dots, n.$

 Using
(\ref{3}), (\ref{3.4}) and the fact that $[e_{i},
e_{j}](H)=0=\nabla_{e_{i}}e_{j}(H)-\nabla_{e_{j}}e_{i}(H)=\omega_{ij}^{1}e_{1}(H)-\omega_{ji}^{1}e_{1}(H),$
for $i\neq j$ and $i, j=2, \dots, n$, we find
\begin{equation}\label{3.9}
\omega_{ij}^{1}=\omega_{ji}^{1}.
\end{equation}

Using (\ref{2.4}), (\ref{4.1}) and $\lambda_{1}=-\frac{nH}{2}$, we
obtain
\begin{equation}\label{4.3}
\sum_{i=2}^{r}{\lambda_{i}}=\frac{n(n-r+3)}{2(n-r+1)}H,
\end{equation}
and $\lambda=\frac{nH}{n-r+1}$.

Therefore, using (\ref{3}) and (\ref{4.3}), we obtain
\begin{equation}\label{4.4}
e_{1}(\lambda_{i})\neq 0,\quad e_{j}(\lambda_{i})=0,
\end{equation}
for $i=1,2,\dots, n$ and $j=2,3,4,\dots, n.$

Now, it can be seen that $\lambda_{1}$ can never be equal to
$\lambda_{i}$ ($i=2, 3, \dots , r$) and $\lambda$.
 Indeed, if $\lambda_{1}= \lambda_{i}$ for some $i$, then from (\ref{3.6}), we
 find that
\begin{equation}\label{3.10}
 e_{1}(\lambda_{j})= (\lambda_{1}-\lambda_{j})\omega_{j1}^{j}=0,
 \quad j=2, 3, \dots ,r,
\end{equation}
which contradicts the first expression of (\ref{4.4}).
Similarly, if $\lambda_1=\lambda$ we get a contradiction.

 Putting
$i\neq 1, j = 1, r+1, \dots , n$ in (\ref{3.6}) and using (\ref{4.4}) and
(\ref{3.5}), we
 find
\begin{equation}\label{4.5}
\omega_{1i}^{1}= \omega_{Ai}^{A}=\omega_{11}^{i}= \omega_{AA}^{i}=
0, \hspace{1 cm} i= 1, 2, A=r+1, \dots , n.
\end{equation}

Putting $k = 1,$ $i\neq j$, and $ i,j=2, 3, \dots ,n$ in
(\ref{3.7}), and using (\ref{3.5}), we
 get
\begin{equation}\label{4.6}
\omega_{ij}^{1}=\omega_{i1}^{j} =\omega_{1i}^{j}
=\omega_{iA}^{1}=\omega_{Ai}^{1} = \omega_{i1}^{A}=
\omega_{A1}^{i}=\omega_{1A}^{i}=\omega_{1i}^{A}= 0, \quad A =r+1, \dots , n.
\end{equation}

Now, putting  $i = 1,2, \dots ,r$, $k = r+1,\dots , n$ and $ j=r+1,
\dots, n$ ($j\ne k$) in (\ref{3.7}), and using (\ref{3.5}), we
 get
\begin{equation}\label{4.7}
\omega_{B1}^{A}=\omega_{AB}^{1} =
\omega_{Bi}^{A}=\omega_{AB}^{i} = 0, \quad
i=1,2, \dots ,r
\end{equation}
where $A\neq B$ and $A, B=r+1,\dots,n.$

Now, evaluating $g( R(e_{1},e_{i})e_{1},e_{i})$, using
(\ref{4.5})$\sim$(\ref{4.7}) and Gauss equation
(\ref{2.5}), we find the following:\\

  $\bullet\hspace{.2 cm}
 \mbox{For}\hspace{.2 cm} X=e_{1}, Y=e_{i}, Z=e_{1}, W=e_{i},$
\begin{equation}\label{3.18}
 e_{1}(\omega_{ii}^{1})- (\omega_{ii}^{1})^{2}= -\frac{nH}{2}
 \lambda_{i}, \quad i=2,3,\dots,n.
\end{equation}

Now, using $\lambda_{1}=-\frac{nH}{2}$,
$\lambda=\frac{nH}{n-r+1}$, and (\ref{3.6}) for $i=1$ and $j=r+1, \dots , n$,
we get
\begin{equation}\label{3.21}
2e_{1}(H)=(n-r+1)H \omega^{1}_{AA}, \quad A=r+1, \dots , n.
\end{equation}

 Now, differentiating (\ref{4.3}) along $e_{1}$ two times
 alternatively
by using (\ref{3.18}) and (\ref{3.21}), we obtain
\begin{equation}\begin{array}{rcl}\label{3.22}
\sum_{i=2}^{r}(2\lambda_{i}+nH)\omega^{1}_{ii}=\frac{n(n-r+3)}{2}H\omega^{1}_{AA},\end{array}
\end{equation}
\begin{equation}\begin{array}{rcl}\label{3.23}
\sum_{i=2}^{r}\Big[2(2\lambda_{i}+nH)(\omega^{1}_{ii})^{2}+
\frac{n(n-r+1)}{2}H\omega^{1}_{ii}\omega^{1}_{AA}
-\frac{nH}{2}\lambda_{i}(2\lambda_{i}+nH)\Big]\\=\frac{n(n-r+3)^{2}}{4}H(\omega^{1}_{AA})^{2}-\frac{n^{3}(n-r+3)}{4(n-r+1)}H^{3},\end{array}
\end{equation}
for $A=r+1, \dots , n$.

 Eliminating
$\lambda_{2}$ from (\ref{3.22}) and (\ref{3.23}) by using
(\ref{4.3}), we obtain
\begin{equation}\begin{array}{rcl}\label{3.24}
\Big[\frac{2n(n-r+2)}{n-r+1}H-2\sum_{i=3}^{r}\lambda_{i}\Big]\omega^{1}_{22}
+\sum_{i=3}^{r}(2\lambda_{i}+nH)\omega^{1}_{ii}=\frac{n(n-r+3)}{2}H\omega^{1}_{AA},\end{array}
\end{equation}
\begin{equation}\begin{array}{rcl}\label{3.25}
\Big[\frac{4n(n-r+2)}{n-r+1}H-4\sum_{i=3}^{r}\lambda_{i}\Big](\omega^{1}_{22})^{2}-nH\Big[\frac{n(n-r+2)}{2(n-r+1)}H-\sum_{i=3}^{r}\lambda_{i}
\Big] \Big[\frac{n(n-r+2)}{n-r+1}H\\-\sum_{i=3}^{r}\lambda_{i}\Big]+
\sum_{i=3}^{r}\Big[2(2\lambda_{i}+nH)(\omega^{1}_{ii})^{2}+
\frac{n(n-r+1)}{2}H\omega^{1}_{ii}\omega^{1}_{AA}\\
-\frac{nH}{2}\lambda_{i}(2\lambda_{i}+nH)\Big]=\frac{n(n-r+3)^{2}}{4}H(\omega^{1}_{AA})^{2}-\frac{n^{3}(n-r+3)}{4(n-r+1)}H^{3},\end{array}
\end{equation}
respectively.

 We consider (\ref{3.24}), (\ref{3.25}) as polynomials of $\omega^{1}_{22}$ with coefficients in polynomial ring
$R_{1}[H,\lambda_{3},\lambda_{4},
\dots, \lambda_{r}, \omega^{1}_{33}, \omega^{1}_{44},
\dots, \omega ^1_{rr}, \omega^{1}_{AA}]$ over real field $\mathbb{R}$. Since
equations (\ref{3.24}), (\ref{3.25}) have a common root
$\omega^{1}_{22}$,  Lemma \ref{lemma2.2} implies that the resultant of
their coefficients is equal to zero, which gives another polynomial
equation defined as
\begin{equation}\label{3.26}
g_{1}(H,\lambda_{3},\lambda_{4},
\dots,\lambda_{r},\omega^{1}_{33},\omega^{1}_{44},
\dots, \omega ^1_{rr}, \omega^{1}_{AA})=0.
\end{equation}

Again differentiating (\ref{3.26}) along $e_{1}$ two times
alternatively and using (\ref{3.18}) and (\ref{3.21}), we obtain
two polynomial equations defined as
\begin{equation}\label{3.27}
g_{2}(H,\lambda_{3},\lambda_{4},
\dots,\lambda_{r},\omega^{1}_{33},\omega^{1}_{44},
\dots, \omega ^1_{rr}, \omega^{1}_{AA})=0,
\end{equation}
\begin{equation}\label{3.28}
g_{3}(H,\lambda_{3},\lambda_{4},
\dots,\lambda_{r},\omega^{1}_{33},\omega^{1}_{44},
\dots, \omega ^1_{rr}, \omega^{1}_{AA})=0.
\end{equation}

 We consider (\ref{3.26}), (\ref{3.27}) and
(\ref{3.26}), (\ref{3.28}) as polynomials of $\lambda_{3}$ with
coefficients in polynomial ring $R_{2}[H,\lambda_{4},\lambda_{5},
\dots,\lambda_{r},\omega^{1}_{33},\omega^{1}_{44},
\dots,\omega^{1}_{AA}]$ over real field $\mathbb{R}$. Also,
Equations (\ref{3.26}), (\ref{3.27}) and (\ref{3.26}), (\ref{3.28})
have a common root $\lambda_{3}$ and Lemma \ref{lemma2.2} implies
that the resultants of their coefficients are equal to zero, which gives
polynomial equations
\begin{equation}\label{3.29}
g_{4}(H,\lambda_{4},\lambda_{5},
\dots,\lambda_{r},\omega^{1}_{33},\omega^{1}_{44},
\dots, \omega ^1_{rr}, \omega^{1}_{AA})=0,
\end{equation}
\begin{equation}\label{3.30}
g_{5}(H,\lambda_{4},\lambda_{5},
\dots,\lambda_{r},\omega^{1}_{33},\omega^{1}_{44},
\dots, \omega ^1_{rr}, \omega^{1}_{AA})=0.
\end{equation}

Similarly, we can eliminate $\omega^{1}_{33}$ from (\ref{3.29}),
(\ref{3.30}) by considering $\omega^{1}_{33}$ as a common root of
(\ref{3.29}), (\ref{3.30}) and by using Lemma \ref{lemma2.2} we obtain
another polynomial equation
\begin{equation}\label{3.31}
g_{6}(H,\lambda_{4}, \lambda_{5},
\dots, \lambda_{r}, \omega^{1}_{44},\omega^{1}_{55},
\dots, \omega ^1_{rr}, \omega^{1}_{AA})=0.
\end{equation}

 Proceeding in the same way, we will be able to
eliminate $\lambda_{4}$, $\omega^{1}_{44}$, $\lambda_{5}$,
$\omega^{1}_{55}$ $\dots$, $\lambda_{r}$, $\omega^{1}_{rr}$,
$\omega^{1}_{AA}$ and obtain a polynomial equation in $H$ with
constant coefficients, which implies that $H$ must be a constant.

\smallskip
\noindent
\textbf{Case B.} $\lambda_{i}= \lambda_{j},$ for some $i,j=2,3,
\dots, r.$\\

For simplicity we will prove it for $i=2, j=3$ and the other cases can be obtained similarly.
Then for $\lambda_{2}= \lambda_{3}$  (\ref{4.3}) reduces to
\begin{equation}\label{3.32}
2\lambda_{3}+\sum_{i=4}^{r}\lambda_{i}=\frac{n(n-r+3)}{2(n-r+1)}H.
\end{equation}

By differentiating (\ref{3.32}) two times along $e_{1}$ and using
(\ref{3.18}) and (\ref{3.21}), we obtain polynomial equations in
$\lambda_{i}$, $\omega^{1}_{ii}$ and $H$. As in the above case, by
using Lemma \ref{lemma2.2} we will be able to find a polynomial
equation in $H$ with constant coefficients which implies that $H$
must be a constant.

\smallskip
\noindent
\textbf{Case C.} $\lambda_{i}= \lambda$ for some $i=1, 2, \dots , r$.

In a similar way with Case B we obtain that $H$ must be constant,
and this concludes the proof of the claim.

Now, since $H$ is constant  it follows from
(\ref{1.1}) that $H{\rm trace}(\mathcal{A}^{2})=0$,
which implies that $H=0$, and
this concludes the proof of  Theorem \ref{main1}. \qed

\medskip
Proposition \ref{main} now follows from the above proof.
Since $M^n$ is a $\delta(r)$-ideal biconservative hypersurface, equation (\ref{1.2}) is satisfied  and we proved  that this implies that $H$ is constant.

\medskip
 \noindent
\textbf{\emph{Acknowledgement}}: \small \emph{The first author is thankful
to Dr. Ram Shankar Gupta for useful discussions and suggestions.  Both authors express their gratitudes to the referee for useful suggestions towards improvement of the paper.}\\

% ------------------------------------------------------------------------

\bibliography{xbib}

% ------------------------------------------------------------------------%
% ------------------------------------------------------------------------%

\medskip
Authors' addresses:
\\
\textbf{Deepika}\\
S N Bose National Centre for Basic Sciences,\\
JD Block,  Sector III, Salt Lake City,\\
Kolkata- 700106, India.\\
\textbf{Email:} sdeep2007@gmail.com, deepika.kumari@bose.res.in\\
\\
{\it Present address:} Romanian Institute of Science and Technology,\\
 Str. Virgil Fulicea nr. 17, 400022,\\
Cluj-Napoca, Romania.\\
\textbf{Email:} deepika@rist.ro\\

\medskip
\noindent
\textbf{Andreas Arvanitoyeorgos}\\
University of Patras,\\
Department of Mathematics,\\
GR-26500 Patras, Greece.\\
\textbf{Email:} arvanito@math.upatras.gr\\

\end{document}